\documentclass[11pt,a4paper]{amsart}
\usepackage{amssymb,amscd,latexsym}
\usepackage[all]{xy}

\theoremstyle{plain}
\newtheorem{theorem}{Theorem}[section]
\newtheorem{lemma}[theorem]{Lemma}
\newtheorem{proposition}[theorem]{Proposition}
\newtheorem{corollary}[theorem]{Corollary}

\theoremstyle{definition}
\newtheorem{definition}[theorem]{Definition}

\numberwithin{equation}{section}

\def\rep{\mathop\mathrm{rep}}
\def\supp{\mathop\mathrm{supp}}
\def\El{\mathop\mathrm{El}}
\def\Dim{\mathop\mathbf{dim}}
\def\Hom{\mathop\mathrm{Hom}\nolimits}
\def\End{\mathop\mathrm{End}}
\def\id{\mathop\mathrm{Id}}
\def\im{\mathop\mathrm{Im}}
\def\Aut{\mathop\mathrm{Aut}}

\def\bA{\mathbf A}	\def\bB{\mathbf B}	\def\bC{\mathbf C}
\def\bS{\mathbf S}	\def\bG{\mathbf G}	\def\fD{\mathbf d}
\def\fX{\mathbf x}	\def\fC{\mathbf c}	\def\sU{\mathsf U}
\def\de{\delta}		\def\De{\Delta}	\def\Th{\Theta}
\def\La{\Lambda}	\def\rQ{\mathrm Q}	\def\dK{\mathfrak K}
\def\Mk{\Bbbk}	\def\mN{\mathbb N}	\def\io{\iota}
\def\rO{\mathrm O}	\def\nW{\mathrm w}	\def\bT{\mathbf T}
\def\kJ{\mathcal J}	\def\tV{\tilde V}	\def\tio{\tilde\io}

\let\tilde=\widetilde  	\let\hat=\widehat
\def\pr{\prec}	\def\pe{\preccurlyeq}
\def\set#1{\{\,#1\,\}}	\def\setsuch#1#2{\{\,#1\,|\,#2\,\}}
\def\sbe{\subseteq}	\def\*{\otimes}	\def\+{\oplus}
\def\bop{\bigoplus}	\def\ol{\overline}	\def\iso{\simeq}
\def\lst#1#2{#1_1,#1_2,\dots,#1_{#2}}
\def\iff{if and only if }		\def\xx{\times}
\def\Arr{\Rightarrow}		\def\={\setminus}
\def\bul{{\bullet}}	\def\0{\emptyset}

\begin{document}

\title[Dimesnions of finite type]{Dimensions of finite type for representations of partially ordered sets}
\author[Y. Drozd and E. Kubichka]{Yuriy A. Drozd and Eugene A. Kubichka}
\address{Department of Mechanics and Mathematics\\ Kyiv Taras Shevchenko University\\ 
 01033 Kyiv\\ Ukraine}
\email{yuriy@drozd.org}
\email{zhenya@gioc.kiev.ua}
\urladdr{drozd.org/$\sim$yuriy}
\subjclass[2000]{16G20, 16G60}
\keywords{Representations of posets, finite type, indecomposable representations}

\begin{abstract}
 We consider the \emph{dimensions of finite type} of representations of a 
 partially ordered set, i.e. such that there is only finitely many isomorphism classes
 of representations of this dimension. We give a criterion for a dimension to be of
 finite type. We also characterize those dimensions of finite type, for which there is an 
 indecomposable representation of this dimension, and show that there can be
 at most one indecomposable representation of any dimension of finite type.
 Moreover, if such a representation exists, it only has scalar endomorphisms.
 These results (Theorem~\ref{main}, page~\pageref{main}) generalize those
 of \cite{k1,d1,mz}.
\end{abstract}

\maketitle

\section{Preliminaries and the Main Theorem}
\label{s1} 

 Let $(\bS,\pr)$ be a finite partially ordered set (\emph{poset}); we denote by $\pr$ the \emph{strict order}
 and by $a\pe b$ the relation ``$a\pr b$ or $a=b$''. We usually suppose that $\bS=\set{1,2,\dots,n}$
 (certainly, not necessary with the natural order) and denote by $\hat\bS =\bS\cup\set0$. (Note that 
 we do not treat $\hat\bS$ as a poset.)  A \emph{representation}
 $V$ of $\bS$ over a field $\Mk$ is an order preserving map of $\bS$ into the set of subspaces of
 a finite dimensional vector space $V(0)$ over $\Mk$. A \emph{morphism} $f:V\to V'$ of such representations
 is a linear mapping $f:V(0)\to V'(0)$ such that $f(V(a))\sbe V'(a)$ for every $a\in\bS$. 
 We denote by $\rep\bS$ the category of such representations (supposing the field $\Mk$ fixed).

 Recall a relation to a bimodule category \cite{d1}. Let $\La=\La_\bS$ be the incidence algebra of the 
 poset $\bS$, i.e. the subalgebra of $\mathrm{Mat}\,(n,\Mk)$ with the basis $\setsuch{e_{ab}}{a\pr b \text{ in }\bS}$.
 Let also $\sU=\sU_\bS$ be the right $\La$-module with the basis $\lst vn$ and the action
 $v_ae_{bc}=\de_{ab}v_c$. We consider $\sU$ as $\Mk$-$\La$-bimodule. Then the category $\El(\sU)$
 of \emph{elements of} $\sU$ (or of \emph{matrices with entries from} $\sU$) is defined. Its objects are
 elements from $\Hom_\La(P,L\*_\Mk\sU)$, where $L$ is a finite dimensional vector space and $P$
 is a finitely generated (right) projective $\La$-module. A \emph{morphism} $\phi:u\to u'$, where
 $u:P\to L\*\sU,\ u':P'\to L'\*\sU$ is, by definition, a pair $\phi_0,\phi_1$, where $\phi_0:L\to L'$
 is a linear map, $\phi_1:P\to P'$ is a $\La$-homomorphism, such that $u'\phi_1=(\phi_0\*1)u$. Given
 an element $u:P\to L\*U$, set $V(0)=L$ and $V(a)=\setsuch{v\in L}{v\*v_a\in\im u}$. We get a
 representation $V=\rho(u)\in\rep\bS$. Obviously, if $\phi=(\phi_0,\phi_1)$ is a morphism $u\to u'$,
 then $\phi_0$ is a morphism $\rho(\phi):\rho(u)\to\rho(u')$. So $\rho$ is a functor $\El(\sU)\to\rep\bS$. It is
 not an equivalence, but one can easily control its defects. Namely, let $\La_a=e_{aa}\La$; they are
 all indecomposable projective $\La$-modules. Consider the so-called \emph{trivial element}
 $T_a$, which is the unique element of $\Hom_\La(\La_a,0\*\sU)$ (it is not zero in the category
 $\El(\sU)$\,).  Later on we shall also use the trivial representation $T_0\in\Hom(0,\Mk\*\sU)$.

  \begin{proposition}\label{11} 
  \begin{enumerate}
\item  The functor $\rho$ is \emph{dense} (i.e. every object from $\rep\bS$ is isomorphic to $\rho(u)$ for some $u$)
 and \emph{full}, i.e. all induced maps $\Hom(u,u')\to\Hom(\rho(u),\rho(u'))$ are surjective.
 \item  $\rho(\phi)=0$ \iff $\phi$ factors through a direct sum $\bop_{a\in \bS}m_aT_a$ of trivial elements.
 In particular, only such direct sums become zero under the functor $\rho$.
 \end{enumerate}
 \end{proposition} 
 \begin{proof} 
 {\em1}. Let $V\in\rep\bS$, $L=V(0)$. Consider the subspace $M=\sum_{a\in\bS}V(a)\*v_a\sbe L\*\sU$.
 It is a $\La$-submodule. Let $P\to M$ be a projective cover of $M$. Considered as a homomorphism
 $P\to L\*\sU$, it defines an element $u\in\El(\sU)$ and it is obvious that $\rho(u)=V$. If $V'=\rho(u')$,
 where $u':P'\to L'\*\sU$, and $f:V(0)=L\to V'(0)=L'$ is a morphism $V\to V'$, then the inclusions $f(V(a))\sbe
 V'(a)$ for all $a\in\bS$ imply that $(f\*1)(\im u)\sbe\im u'$. Hence, there is a homomorphism $g:P\to P'$
 with $(f\*1)u=u'g$, which gives a morphism $\phi=(f,g)$ such that $\rho(\phi)=f$.

 {\em2}. If $\rho(\phi)=0$, then $\phi=(0,\phi_1)$, so $u'\phi_1=0$ and  it decomposes as
 $$ 
   \begin{CD}
     P @>u>> L\*\sU \\
	@|  @VVV \\
     P @>>> 0\*\sU \\
       @V\phi_1VV @VVV\\
     P'@>>u'> L'\*\sU.
 \end{CD}  
 $$ 
 Obviously, the second row of this diagram splits in $\El(\sU)$ into a direct sum of trivial representations.
 \end{proof} 
 
 Note that $\Hom_\La(P_a,\sU)\iso \sU e_{aa}=\langle v_a\rangle$. Therefore, a homomorphism $d_aP_a\to L\*\sU$
 can be identified with a matrix $M(a)$ of size $d_0\xx d_a$, where $d_0=\dim L$. Since every projective
 $\La$-module $P$ decomposes uniquely as $\bop_{a\in\bS}d_a\La_a$, and 
 $$ 
  \Hom_\La(\La_b,\La_a)= \begin{cases}
  \Mk	&\text{ if }\, a\pr b,\\
   0	&\text{ otherwise},
 \end{cases}  
 $$ 
 it gives the original ``matrix'' definition of \cite{nr}. Namely, $u$ is presented as a block matrix
  \begin{equation}\label{e1} 
 M=  \begin{array}{|c|c|c|c|c}
  \cline{1-4} &&&&\\[-1ex] M(1) & M(2) & \dots &M(n)& \\[-1ex] &&&&,\\  \cline{1-4}
 \end{array} 
 \end{equation}
 where $M(a)$ is of size $d_0\xx d_a$. For two matrices of this shape, $M$ and $M'$, a morphism
 $\Phi$ is given by a set of matrices $\{\,\Phi(a)\,|\,a\in\hat\bS\,\}\cup\setsuch{\Phi(ba)}{b\pr a\text{ in }\bS}$
 such that, for every $a\in\bS$, 
  \begin{equation}\label{e2} 
 \Phi(0)M(a)=M'(a)\Phi(a)+\sum_{b\pr a}M'(b)\Phi(ba).  
 \end{equation}

 In some respect, this bimodule (or matrix) interpretation has certain
 advantage, and we shall permanently use it. Especially, it gives rise to a quadratic form useful in many questions.
 
  \begin{definition}\label{12} 
 \begin{enumerate}
\item The \emph{dimension} (or \emph{vector dimension}) of an element $u\in \Hom_\La(P,L\*\sU)$, or of the
 corresponding representation of $\bS$, is the function $\fD=\Dim u:\hat\bS \to\mN$ such that
 $\fD(0)=\dim L$ and $P\iso\bop_{a\in\bS}\fD(a)\La_a$. We denote by $\El_\fD(\sU)$ the set of all
 elements of dimension $\fD$ and by $\rep_\fD(\bS)$ the set of the corresponding representations.

 If $u$ arises as above from a representation $V\in\rep\bS$, then $\fD(0)=\dim V(0)$ and
 $\fD(a)=\dim\big(V(a)/\sum_{b\pr a}V(b)\big)$ for $a\in\bS$.

 \item  The \emph{support} of a dimension $\fD:\hat\bS\to\mN$ is the subset
 $\supp\fD=\setsuch{a\in\bS}{\fD(a)\ne0}$. The dimension $\fD$, as well as the elements from $\El_\fD(\bS)$
 and the corresponding representations, is called \emph{sincere} if $\supp\fD=\hat\bS$.

 If a dimension $\fD$ is not sincere, the representations of this dimension can (and usually will) be treated
 as representations of a smaller poset, namely its support.

 \item    The \emph{quadratic form $\rQ_\bS$ associated to a poset} $\bS$ is, by definition, the quadratic form
 $$ 
  \rQ_\bS(x_0,x_1,\dots,x_n) = \sum_{a\in\hat\bS } x_a^2 +\sum_{\substack{a,b\in\bS\\ a\pr b}}
	 x_ax_b -\sum_{a\in\bS}x_0x_a.
 $$ 
\end{enumerate}
 \end{definition}
 
 Note that if $\fD:\hat\bS \to\mN$, then the negative part of $\rQ_\bS(\fD)$ is just the dimension of the vector
 space $\El_\fD(\sU)=\Hom_\La(P,L\*\sU)$ of all elements of dimension $\fD$, while the positive part is the
 dimension of the algebraic group $\bG_\fD=\Aut L\xx\Aut P$ acting on $\El_\fD(\sU)$ so that its orbits are the
 isomorphism classes of elements. From here the following result is evident.

 \begin{proposition}\label{13} 
 \begin{enumerate}
 \item   If a dimension $\fD:\hat\bS \to\mN$ is \emph{of finite type}, i.e. 
 there are only finitely many isomorphism classes in $\El_\fD(\sU)$, then $\rQ_\bS(\fD')>0$ for each
 dimension $\fD'\le\fD$, i.e. such that $\fD'(a)\le\fD(a)$ for all $a\in\hat\bS $.
 \item   Especially, if $\bS$ is \emph{representation finite}, i.e. has only finitely many nonisomorphic
 indecomposable representations, the quadratic form $\rQ_\bS$ is \emph{weakly positive}, i.e.
 $\rQ_\bS(\fX)>0$ for every nonzero vector $\fX$ with non-negative entries.
\end{enumerate}
 \end{proposition}
 
 In \cite{k1,d1} the converse was proved, giving a criterion for $\bS$ to be representation finite.
 We recall this result. A poset $\bS$ is called \emph{primitive} if it is a disjoint unit of several chains such
 that the elements of different chains are noncomparable. We denote such a poset by $(n_1,n_2,\dots,n_s)$,
 where $n_i$ are the lengths of the chains. We also denote by $\dK$ the poset
 $\set{a_1,a_2,b_1,b_2,c_1,c_2,c_3,c_4}$, where the order $\pr$ is defined as follows:
 $ a_2\pr a  _1,\  b_2\pr b_1,\ b_2\pr a_1,\ c_1\pr c_2\pr c_3\pr c_4$. The posets $(1,1,1,1)$,\,
 $(2,2,2)$,\,$(1,3,3)$,\,$(1,2,5)$ and $\dK$ are called \emph{critical}. 

  \begin{theorem}\label{14} 
 \begin{enumerate}
 \item    The following conditions are equivalent:
	\begin{enumerate}
	 \item   $\bS$ is representation finite.
	 \item  $\rQ_\bS$ is weakly positive.
	 \item  $\bS$ contains no critical subset.
	\end{enumerate}
  \item  Let $\bS$ is representation finite, $\fD:\hat\bS \to\mN$.
 The following conditions are equivalent:
	\begin{enumerate}
	 \item  There is an indecomposable element $u\in\El_\fD(\sU)$.
	 \item  $\fD$ is a \emph{root} of the form $\rQ_\bS$, i.e. $\rQ_\bS(\fD)=1$.
	\end{enumerate}
 Moreover, if the latter condition holds, there is a unique indecomposable element
 $u\in\El_\fD(\sU)$, $\End u=\Mk$ and the orbit of $u$ is open in the space $\El_\fD(\sU)$
 (in the Zariski topology).
 \end{enumerate}
 \end{theorem}

 We shall generalize this result using the following notions.

 \begin{definition}\label{15} 
 Let $\fD:\hat\bS \to\mN$.
 \begin{enumerate}
 \item  The dimension $\fD$ is called \emph{critical}, if its support $\bC$ is a critical subset,
 $\rQ_\bC(\fD)=0$ and the values $\setsuch{\fD(a)}{a\in\hat\bS }$ are coprime (equivalently,
 at least one of these values equals $1$).

 Table~1 below presents all critical dimensions (there are 5 of them, denoted by $\fC_i,\ 1\le i\le5$).
 In every picture from this table the bullets show the elements $a\in\bC$;
 the numbers nearby are the values $\fC_i(a)$. The relations $a\pr b$ are shown by the edges going from
 $a$ downstairs to $b$ upstairs. The number in a circle above denotes the dimension $\fC_i(0)$.
\end{enumerate}
 \end{definition}
 
 \begin{table}[!ht]
  \caption{\bf Critical dimensions}\label{t1} 
 \def\sm#1{{\,\scriptsize#1\,}}
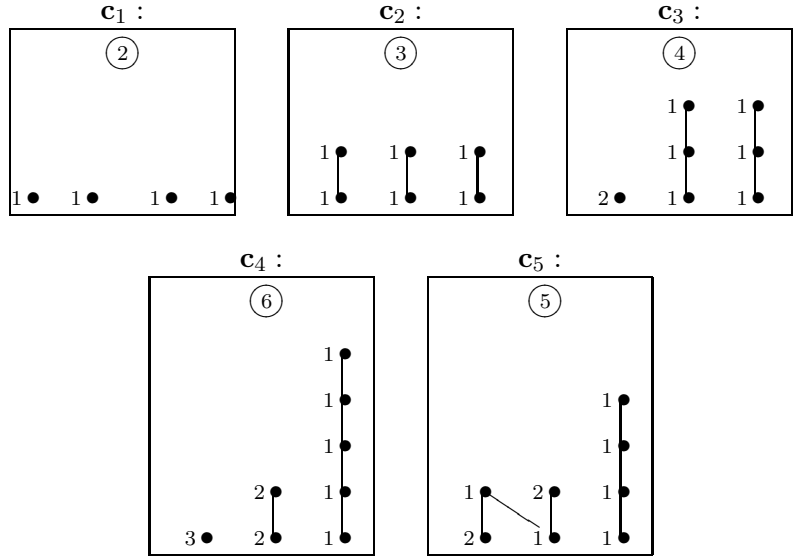
 \begin{center}
 \begin{tabular}{ccccc} 
 $\fC_1:$ && $\fC_2:$ && $\fC_3:$
\\
 \framebox[3cm]{$$\xymatrix@C=1ex@R=1em{&&& *+[o][F-]\txt{\sm2} &&&\\ \\ \\
 *=0+\txt{\sm1$\bul$}&&*=0+\txt{\sm1$\bul$}&&*=0+\txt{\sm1$\bul$}&&*=0+\txt{\sm1$\bul$}
 }$$}&& 
 \framebox[3cm]{$$\xymatrix@C=1ex@R=1em{&& *+[o][F-]\txt{\sm3} \\ \\
	*=0+\txt{\sm1$\bul$} && *=0+\txt{\sm1$\bul$} && *=0+\txt{\sm1$\bul$} \\
	*=0+\txt{\sm1$\bul$}\ar@<-.5ex>@{-}[u] && *=0+\txt{\sm1$\bul$} \ar@<-.5ex>@{-}[u]
	&& *=0+\txt{\sm1$\bul$} \ar@<-.6ex>@{-}[u]
}$$} && 
 \framebox[3cm]{$$\xymatrix@C=1ex@R=1em{ && *+[o][F-]\txt{\sm4} \\
	&& *=0+\txt{\sm1$\bul$} && *=0+\txt{\sm1$\bul$} \\
	&& *=0+\txt{\sm1$\bul$}\ar@<-.5ex>@{-}[u] && *=0+\txt{\sm1$\bul$}\ar@<-.5ex>@{-}[u] \\
	*=0+\txt{\sm2$\bul$}&& *=0+\txt{\sm1$\bul$}\ar@<-.5ex>@{-}[u] && *=0+\txt{\sm1$\bul$}\ar@<-.5ex>@{-}[u] 	
 }$$} \\
  \end{tabular}

\bigskip
\begin{tabular}{ccc}
 $\fC_4:$ &&$\fC_5:$\\
 \framebox[3cm]{$$ \xymatrix@C=1ex@R=1em{ &&*+[o][F-]\txt{\sm6} \\  
	&& && *=0+\txt{\sm1$\bul$} \\
	&& && *=0+\txt{\sm1$\bul$}\ar@<-.5ex>@{-}[u]  \\
	&& && *=0+\txt{\sm1$\bul$}\ar@<-.5ex>@{-}[u]  \\
	&& *=0+\txt{\sm2$\bul$}  && *=0+\txt{\sm1$\bul$}\ar@<-.5ex>@{-}[u]  \\
	*=0+\txt{\sm3$\bul$}	&& *=0+\txt{\sm2$\bul$}\ar@<-.5ex>@{-}[u]   && *=0+\txt{\sm1$\bul$}
	 \ar@<-.5ex>@{-}[u]  }$$} &\qquad&
 \framebox[3cm]{$$
 \xymatrix@C=1ex@R=1em{ &&*+[o][F-]\txt{\sm5} \\ \\
	&& && *=0+\txt{\sm1$\bul$} \\
	&& && *=0+\txt{\sm1$\bul$}\ar@<-.5ex>@{-}[u]  \\
	*=0+\txt{\sm1$\bul$}	&& *=0+\txt{\sm2$\bul$}    && *=0+\txt{\sm1$\bul$} \ar@<-.5ex>@{-}[u]\\
	*=0+\txt{\sm2$\bul$}\ar@<-.5ex>@{-}[u]	&& *=0+\txt{\sm1$\bul$}\ar@<-.5ex>@{-}[u] 
	\ar@<-.5ex>@{-}[ull]  && *=0+\txt{\sm1$\bul$}\ar@<-.5ex>@{-}[u]  
	} $$}\\ 
 \end{tabular} 
\end{center}	
 \end{table}	
 
   \begin{theorem}[Main Theorem]\label{main}
  \begin{enumerate}
 \item  The following conditions for a dimension $\fD:\hat\bS \to\mN$ are equivalent:
	\begin{enumerate}
	 \item  $\fD$ is a dimension of finite type.
	 \item  $\rQ_\bS(\fD')>0$ for every nonzero dimension $\fD'\le\fD$.
	 \item  There is no critical dimension $\fC\le\fD$.
	\end{enumerate}
 \item  If a dimension $\fD$ is of finite type, the following conditions are equivalent:
	\begin{enumerate}
	\item  There is an indecomposable element $u\in\El_\fD(\sU)$.
	 \item  $\rQ_\bS(\fD)=1$. 
	\end{enumerate}
 Moreover, if the latter condition holds, there is a unique indecomposable element $u\in\El_\fD(\sU)$,
   $\End u=\Mk$, and the orbit of $u$ is open and dense in $\El_\fD(\sU)$.
\end{enumerate}
 \end{theorem}

 Note that all claims about indecomposable elements of $\El_\fD(\sU)$ obviously remain valid for
 indecomposable \emph{representations} from $\rep_\fD(\bS)$, with the exception of trivial dimensions,
 which are nonzero on a unique element $a\in\bS$.

 For primitive posets this theorem was deduced in \cite{mz} from the results of Kac \cite{kac} about 
 the representations of quivers. Unfortunately, this approach cannot be applied in general 
 case. That is why we have to return to the original technique of \emph{derivations} (or
 \emph{differentiation}) from \cite{nr}. It will be considered in the next section.

 \section{Derivations and integration}
 \label{s2} 

 For calculation of representations there is an effective algorithm of \emph{derivations} (or differentiation)
 elaborated in \cite{nr,k2}. We recall it; moreover, we show that it can be considered as an equivalence of
 certain categories. For every element $a\in\bS$, denote by $\De(a)=\setsuch{b\in\bS}{b\pe a}$ the
 \emph{lower cone} of $a$, $\De'(a)=\De(a)\=\set{a}$, and $\Th(a)$ the set of elements
 noncomparable with $a$. Let also $\nW(\bS)$ be the \emph{width} of $\bS$, i.e. the maximal number
 of pairwise noncomparable elements from $\bS$.
 
  \begin{definition}\label{21} 
  Suppose that $a$ is a maximal element of $\bS$.
     Let $\Pi(a)$ be the set of all pairs $\set{b,c}$ such that $b,c\in\Th(a)$ and are noncomparable in $\bS$.
  Set $\tilde\bS^a=\bS\cup \Pi(a)$ and define a partial order $\pe$
 on $\tilde\bS^a$ setting $B\pe C$ in $\tilde\bS^a$ \iff for each element $b\in B$ there is an element $c\in C$
 such that $b\pe c$ in $\bS$ (we identify elements of $\bS$ with one-element sets). We also set
 $\bS^a=\tilde\bS^a\=\set a$ and call the poset $\bS^a$ the \emph{derivative of $\bS$ with respect to $a$}.

 \smallskip
 For instance, $b\pe\set{c,d}$ means that either $b\pe c$ or $b\pe d$; $\set{b,c}\pe d$ means that both
 $b\pe d$ and $c\pe d$, etc.

 \smallskip
 We fix, for every pair $p\in\Pi(a)$, one element $p'\in p$, and denote by $p''$ the other element of
 $p$.
 \end{definition}

 We also use the following notations.
 \begin{itemize}
 \item  For every element $a\in\bS$ denote by $E_a$ the representation of $\bS$ such that
 $E_a(0)=\Mk$, $E_a(b)=\Mk$ if $a\pe b$ and $E_a(b)=0$ otherwise.
 \item  For every pair of noncomparable elements $p=\set{a,b}$ of $\bS$, denote by $E_p$ the
 representation of $\bS$ such that $E_p(0)=\Mk$, $E_p(c)=\Mk$ if $a\pe c$ or $b\pe c$
 and $E_p(c)=0$ otherwise.
\end{itemize}

 We use the same notations for the objects of $\El(\sU)$ corresponding to these representations. In the
 matrix form, $E_a(a)=(1)$, $E_a(b)=\0$ if $b\ne a$; $E_p(a)=E_p(b)=(1)$, $E_p(c)=\0$ if
 $c\ne a,\,c\ne b$. 

\medskip
 If $V$ is a representation of $\bS$, define the \emph{derived representation} $D_aV$ of $\bS^a$ as
 follows:
 \begin{itemize}
\item 	$D_aV(0)=V(a)$;
 \item  $D_aV(b)=V(b)\cap V(a)$ for $b\in\bS\=\set a$;
 \item  $D_aV(p)=(V(p')\cap V(p''))\cap V(a)$ for $p\in \Pi(a)$.
\end{itemize}
 Obviously, every morphism $f:V\to W$ induces a morphism $D_af:D_aV\to D_aW$. So we obtain a
 functor $D_a:\rep(\bS)\to\rep(\bS^a)$.

 On the contrary, let $V$ be a representation of $\bS^a$. For every $p\in \Pi(a)$, let $\tV(p)=V(p)/(V(p')+V(p''))$
 and $\pi_p:V(p)\to\tV(p)$ be the natural surjection. We can choose sections $\io_p:\tV(p)\to V(p)$ such that
 $\pi_p\io_p=\id$ and $\io_p|_{V(q)}=\io_q$ if $p,q\in \Pi(a),\ q\pr p$.
 Set $\tV(0)=\bop_{p\in \Pi(a)}\tV(p)$ and define, for $b\in\Th(a)$, a map $\tio_b:\tV(0)\to V(0)\+\tV(0)$
 by the rule
 $$ 
 \tio_b(v) = \begin{cases}
  (0,0) &\text{ if } b\notin p,\\
  (0,v) &\text{ if } b=p',\\
 (\io_p(v),v) &\text{ if } b=p'',
 \end{cases}
 $$ 
 where $v\in\tV(p)$.
 We construct the \emph{integrated representation} $\int_a^\io V$ as follows:
 \begin{itemize}
\item 	$\int_a^\io V(0)=V(0)\+\tV(0)$;
 \item  $\int_a^\io V(a)=V(0)$;
 \item  $\int_a^\io V(b)=V(b)\,$ for $b\in\De(a)$;
 \item  $\int_a^\io V(b)=V(b)+\im\tio_b\,$ for $b\in\Th(a)$.
\end{itemize}
 We have included the choice $\io=\set{\io_p}$ of sections $\io_p:\tV(p)\to V(p)$ into this notation. 
 Nevertheless, if $\io'=\set{\io'_p}$ is another choice of such sections, $\im(\io'_p-\io_p)\in V(p')+V(p'')$
 for each $p\in \Pi(a)$. Thus we can find maps $\de_p:\tV(p)\to V(p')$ such that $\im(\io'_p-\io_p-\de_p)
 \sbe V(p'')$. Moreover, we can again suppose that $\de_p|_{\tV(q)}=\de_q$ if $p,q\in \Pi(a),\ q\pr p$.
 It defines a map $\de:\tV(0)\to V(0)$ such that the map $V(0)\+\tV(0)\to V(0)\+\tV(0)$
 given by the matrix
 $$ 
   \begin{pmatrix}
  \id & \de \\ 0 &\id
 \end{pmatrix}  
 $$ 
 is indeed a morphism (hence, an isomorphism) $\int^\io_aV\to\int^{\io'}_aV$. So we can use the notation
 $\int_aV$ without mentioning $\io$. Note that we have only defined the operation $\int_a$ on
 representations, not on their morphisms, so it is not a functor. Nevertheless, Proposition~\ref{23} below
 shows that it can be considered as a functor from $\rep\bS^a$ to a factorcategory of $\rep\bS$.

 This integration is easier in the matrix language. Namely, let a set of matrices $\setsuch{M(x)}{x\in\bS^a}$
 define an object $u\in\El(\sU_{\bS^a})$, like in \eqref{e1}, and $\fD=\Dim u$. Let also $\De(a)=\set{a=\lst ak}$.
 We choose a matrix $M(a)$ with $\fD(0)$ rows so that its columns are linear independent and
 $$ 
  \mathrm{rank}\  \begin{array}
  {|c|c|c|c|}
  \hline &&&\\[-1ex] M(a_1) & M(a_2) & \dots &M(a_n) \\[-1ex] &&& \\ \hline
 \end{array}  =\fD(0).
 $$ 
 We denote by $\fD(a)$ the number of columns of $M(a)$. Define $\fD^*:\hat\bS \to\mN$ as follows:
 $$ 
  \fD^*(b)= \begin{cases}
   \fD(0)+\sum_{p\in\Pi(a)}\fD(p) &\text{if }\, b=0,\\
   \fD(b) &\text{if }\, b\in\De(a),\\
   \fD(b)+\sum_{\substack{ p\in\Pi(a)\\ b\in p\ }} \fD(p) &\text{if }\,b\in\Th(a).
 \end{cases}  
 $$ 
 The integrated element $\int_au$ is of dimension $\fD^*$ and is given by the set of matrices $M^*(b)$
 defined as follows. We consider the element $z\in\El(\sU)$, which is the direct sum
 $$
 z=\big(\mspace{-9mu}\bop_{b\in\Th(a)}\fD(b)E_b\big) \+\big(\mspace{-9mu}\bop_{p\in\Pi(a)}E_p\big).
 $$
 In the block matrix $Z$ defining this element only blocks $Z(b),\ b\in \Th(a)$, are nonzero; let 
 \begin{align*} 
    Z(b)&=  \begin{array}{|c|c|c|c|c}
  \cline{1-4} &&&&\\[-1ex] Z_b(0) & Z_b(p_1) & \dots &Z_b(p_s)& \\[-1ex] &&&&,\\  \cline{1-4} 
 \end{array} \\ \\
   Y(b)&=  \begin{array}{|c|c|c|c|c}
  \cline{1-4} &&&&\\[-1ex] M(b) & M_b(p_1) & \dots &M_b(p_s)& \\[-1ex] &&&&,\\  \cline{1-4} 
 \end{array}
 \end{align*}  
 where
 \begin{itemize}
\item  $\lst ps$ are all pairs from $\Pi(a)$ containing $b$; $Z_b(p_i)$ denotes the part of $Z(b)$
 corresponding to the direct summand $E_{p_i}$ of $Z$, and $Z_b(0)$ is the part of $Z(b)$
 corresponding to the direct summand $E_b$ (it is the zero matrix with $\fD(b)$ columns);
 \item  the vertical stripes of the matrix $Y(b)$ are of the same size as the corresponding stripes
 of the matrix $Z(b)$;
 \item  $M_b(p)=M(p)$ if $b=p''$, and $M_b(p)=0$ if $b=p'$.
\end{itemize}
We also set $Z(b)=0$ if $b\in\De(a)$. Then
 $$ 
  M^*(b)= \begin{array}{|c|l}
  \cline{1-1} \\[-2ex] Y(b)& \\ \cline{1-1} \\[-2ex] Z(b)&.\\\cline{1-1}
 \end{array}   
 $$ 

 Since $\nW(\Th(a))\le2$, every object of $\rep\bS$ with support in $\Th(a)$ is a direct
 sum of the trivial representation $T_0$, the representations $E_b$ and $E_p$. Set
 $\rO(a)=\setsuch{T_0,E_b,E_p}{b\in\Th(a),\,p\in\Pi(a)}$. They are all indecomposable
 representations $V$ such that $D_aV=0$. Straightforward matrix calculations
 immediately imply the following result (cf. also \cite{k2} and, for paragraphs \emph{3}
 and \emph{4}, the proof of Lemma~\ref{44} below).

 \begin{proposition}\label{23} 
  \begin{enumerate}
\item  If $V\in\rep(\bS^a)$, then $D_a\int_aV\iso V$.
 \item If $V\in\rep(\bS)$, then $\int_aD_aV\iso V $ \iff $V$ has no direct summands from
 $\rO(a)$.
 \item  For every morphism $\phi:V\to W$ of representations of $\bS^a$, there is a morphism
 $f:\int_aV\to\int_aW$ such that $\phi=D_a f$. If, moreover, $\phi$ is an isomorphism, so is $f$. 
 \item  The operations $D_a$ and $\int_a$ induce an equivalence between the categories
 $\rep\bS/\kJ_a$ and $\rep\bS^a$, where $\kJ_a$ is the ideal generated by the identity morphisms
 of all representations from $\rO(a)$.
 \end{enumerate}
 \end{proposition}
 
 We shall call a dimension $\fD':\hat\bS^a\to\mN$ \emph{subordinate} to a dimension $\fD:\hat\bS\to\mN$
 if $\Dim\int_aV=\fD$ for some representation $V\in\rep_{\fD'}(\bS^a)$. Obviously, for every
 dimension $\fD:\hat\bS\to\mN$ there is only a finite set of subordinate dimensions $\fD':\hat\bS^a\to\mN$.
 Proposition~\ref{23} immediately implies the following corollary.

 \begin{corollary}\label{24} 
  If  a dimension $\fD'$ is subordinate to a dimension $\fD$, which is of finite type, then
 $\fD'$ is of finite type as well.
 \end{corollary}

  \section{Dimensions of finite type}
 \label{s3} 

 In this section we shall prove paragraph \emph{1} of the Main Theorem \ref{main}.
 In fact, \emph{1(a)}$\Arr$\emph{1(b)} is the claim of Proposition~\ref{13}.1, and 
 \emph{1(b)}$\Arr$\emph{1(c)} is obvious. So we only have to prove
 that \emph{1(c)}$\Arr$\emph{1(a)}.  

  \begin{definition}\label{31} 
  We call a representation $V$ \emph{quite sincere} if it is indecomposable and the
 following conditions hold:
 \begin{itemize}
 \item   $V(a)\ne V(0)$\, for every $a\in\bS$.
 \item   $V(a)\ne\sum_{b\pr a}V(b)$\, for every $a\in\bS$.
\end{itemize}
 In particular, since $V$ is indecomposable, $\sum_{a\in\bS}V(a)=V(0)$.
 If there is a quite sincere representation of dimension $\fD$, we call this dimension
 \emph{quite sincere} as well.

 Obvious necessary conditions for a dimension $\fD$ to be quite sincere are:
\begin{itemize}
\item  $\fD(a)\ne0$\, for every $a\in\bS$;
 \item  $\sum_{i=1}^k\fD(a_i)<\fD(0)$\, for every chain $a_1\pr a_2\pr \dots \pr a_k$ from $\bS$.
\end{itemize}
 Note also that if there is a quite sincere dimension of representations of $\bS$, $\bS$ must have
 at least $2$ maximal elements.
 \end{definition}

  We shall deduce the implication \emph{1(c)}$\Arr$\emph{1(a)} from the following result. 

  \begin{lemma}\label{32} 
 Suppose that $\nW(\bS)\le3$ and a quite sincere dimension $\fD$ satisfies condition 1(c)
 of Theorem~\ref{main}. There is a maximal element $a\in\bS$ such that every dimension
 $\fD'$ of representations of the derived poset $\bS^a$, which is subordinate to $\fD$,
 satisfies this condition too.
 \end{lemma} 
 \begin{proof} 

 Note, first of all, that $\fD'|_{\bS\=\set a}\le\fD|_{\bS\=\set a}$. Hence, if $\fD'\ge\fC$ for a
 critical dimension $\fC$, the support of this $\fC$ must contain at least one element from $\Pi(a)$:
 otherwise also $\fD\ge\fC$. If there is a maximal element $a\in\bS$ such that $\nW(\Th(a))\le1$,
 then $\bS^a=\bS\=\set a$, so there is nothing to prove. Hence, we may suppose that
 $\nW(\Th(a))=2$ for every maximal element $a\in\bS$. We show that then $\bS$ must have $3$
 maximal elements. Indeed, suppose that  $\bS$ has only $2$ maximal elements, $a$ and $b$.
 If both $a$ and $b$ can be included in non-comparable triples, respectively, $\set{a,a',a''}$
 and $\set{b,b',b''}$, the quadruple $\set{a,a'',b',b''}$ is non-comparable too, which contradicts
 the condition. So, for one of these elements, say for $a$, $\nW(\Th(a))=1$, the case already excluded.

\medskip
 We also recall the following (rather easy) lemma from \cite{nr}.

  \begin{lemma}[\cite{nr}]\label{semidec}  
  Suppose that $\bS=\bS_1\sqcup\bS_2\sqcup\bS_3$, where $\bS_3$ is a chain (maybe empty),
 $b\pr a$ for every $a\in\bS_1,\,b\in\bS_2$ and $\bS_1\ne\0,\,\bS_2\ne\0$. If $V\in\rep\bS$
 is indecomposable and $\fD=\Dim V$, then either $\fD|_{\bS_1}=0$ or $\fD|_{\bS_2}=0$.
 Especially, if neither $\bS_1$ nor $\bS_2$ are empty, $\bS$ has no sincere indecomposable
 representations.

 \em In this case the poset $\bS$ is called \emph{semidecomposable}. Thus in what follows we
 may suppose that $\bS$ is not semidecomposable. 
 \end{lemma}

 Let the maximal elements of $\bS$ be $a,b,c$. Using the Dilworth theorem
 \cite[Theorem 10.2.3]{ore}, we consider $\bS$ as a union of three chains 
 \begin{align*} 
 \bA&=\set{a=a_1\succ a_2\succ\dots\succ a_r},\\
 \bB&=\set{b=b_1\succ b_2\succ\dots\succ b_s},\\
 \bC&=\set{c=c_1\succ c_2\succ\dots\succ c_t},
 \end{align*} 
 Since $\fD$ is quite sincere, $\bS$ contains no primitive subset of type $(2,2,2)$. 
 Consider the \emph{top} of $\bS$, i.e. the maximal primitive subset $\bT\sbe\bS$
 containing $\set{a,b,c}$. Then $\bT\iso(1,m,n)$, namely, we may suppose that
 $\bT=\set{a,\,\lst bm,\,\lst cn}\ (m\le s,\,m\le n\le t) $. If $m=1$, the derived poset $\bS^c$ only
 has one new point $p=(a,b)$, such that $p\succ a_i$ and $p\succ b_j$ for all $i,j$.
 So $p$ cannot occur in any critical subset. If $m\ge3$ and $n\ge3$, then $\fD(a)=1$.
 Hence, in any subordinate dimension $\fD'$ of the poset $\bS^c$, which consists of $\bS\=\set c$
 and the points $p_j=(a,b_j)\ (1\le j\le m)$, only one of the points $p_j$ can occur with
 $\fD'(p_j)=1$, and then $\fD'(a)=0$. Thus, replacing this $p_j$ by $a$, we get the same
 dimension for $\bS$, so $\fD'\ge\fC_i$ is impossible. 
 Therefore, we may suppose that $m=2,\,n\ge2$.

 We distinguish the following cases.

\medskip
 {\bf Case 1}. Either $r=1$ or $a_2\pr c$.

 Then $\bS^c=(\bS\=\set c)\cup\set{(a,b),(a,b_2)}$ and $(a,b)$ cannot occur in any critical subset
 of $\bS^c$:
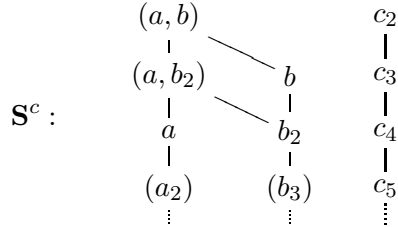
 $$ 
 \bS^c:\qquad \vcenter{
 \xymatrix@=1ex{ {(a,b)}\ar@{-}[d] \ar@{-}[drr] && && c_2 \ar@{-}[d]\\
	(a,b_2) \ar@{-}[d] \ar@{-}[drr] && b \ar@{-}[d] && c_3 \ar@{-}[d] \\
  	   a	\ar@{-}[d] && b_2 \ar@{-}[d] && c_4 \ar@{-}[d] \\
	  (a_2)\ar@{.}[d] && (b_3)\ar@{.}[d] && c_5\ar@{.}[d]	\\
	 && &&  } }
 $$ 
 ($a_2$ and $b_3$ are definitely not in the top of $\bS$). 
 Set $p=(a,b_2)$ and suppose that $\fD'\ge\fC_i$ for some $i$, $\bS_0=\supp\fC_i$.
 It is easy to see that $\bS_0\ne(2,2,2)$. If it is $(1,3,3)$, it can only be $\set{b,\,p\succ a\succ a_2,
 \,c_j\succ c_{j+1}\succ c_{j+2}}$. Then $\set{a\succ a_2,\,b\succ b_2,\, c_j\succ c_{j+1}}$ is a subset of $\bS$
 of type $(2,2,2)$, which is impossible. Suppose that $\bS_0$ is $(1,2,5)$, so $\fD'=\fC_4$. It can only
 occur as $\set{b,\,p\succ a,\,c_j\succ c_{j+1}\succ\dots\succ c_{j+4}}$ with $\fD'(p)\ge2,\,\fD'(a)\ge2,\,
 \fD'(b)\ge3$. Then $\set{a,\,b\succ b_2,\,c_j\succ\dots\succ c_{j+4}}$ is also of type $(1,2,5)$ and
 $\fD(a)\ge4,\, \fD(b)\ge3,\,\fD(b_2)\ge2$, so $\fD\ge\fC_4$, which is impossible. Just in the same
 way, if $\bS_0\iso\dK$, it can only be $\set{a\pr p\succ b_2\pr b,\,c_j\succ\dots\succ c_{j+3}}$ with
 $\fD'(a)\ge2,\,\fD(b)\ge2$. Then $\fD(a)\ge3,\,\fD(b_2)\ge2$, hence $\fD\ge \fC_4$ with 
 $\supp\fC_4=\set{a,\,b\succ b_2,\,c\succ c_j\succ\dots\succ c_{j+3}}$. It is also impossible, which
 accomplishes the consideration of Case 1.

\medskip
 {\bf Case 2}. $a_2\pr b_2$. 

 Then, if $b_3\pr a$, $\bS$ is semidecomposable with $\bS_1=\set{a,b,b_2},\ \bS_2=\setsuch{a_i,b_j}{i>1,j>2}$,
 both nonempty, so there are no quite sincere dimensions at all. Hence either $s=2$ or $b_3\le c$. In both cases
 $\bS^c$ is as in Case 1 and analogous considerations prove the lemma.

\medskip
 {\bf Case 3}. $a_2\pr b,\ a_2\nprec b_2,\ a_2\nprec c$.

 Then the new elements in $\bS^b$ are
 $p_i=(a,c_i),\ (1\le i\le n)$ and $p_1$ cannot occur in any critical subset:
 $$ 
  \bS^b:\qquad \vcenter
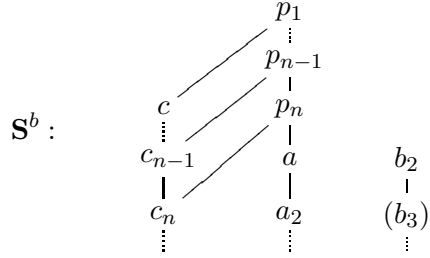
{
 \xymatrix@=1ex {	&& p_1 \ar@{.}[d] \ar@{-}[ddll] && \\
	 	 &&{\ p_{n-1}} \ar@{-}[d] \ar@{-}[ddll] &&\\
		c \ar@{.}[d] && p_n \ar@{-}[d]\ar@{-}[ddll] && \\
		{\ c_{n-1}} \ar@{-}[d] && a \ar@{-}[d] && b_2 \ar@{-}[d] \\
		c_n \ar@{.}[d]  && a_2\ar@{.}[d] && (b_3)\ar@{.}[d]\\
		&&&&	
}}
 $$ 
 ($a_2$ is in the top of $\bS^b$, while $b_3$ is not). Note that either $b_3\pr a$ or $b_3\pr c_{n-1}$:
 otherwise $\bS$ contains a subset $(2,2,2)$. It implies that $b_3$ cannot occur in a critical
 subset $\bS_0\sbe\bS^b$ containing a new element $p_i$. 
 Hence, $\bS_0\ne(2,2,2)$ and $\bS_0\ne\dK$.
 Suppose that $\fD'\ge\fC_i$ with $\supp\fC_i=\bS_0$.
 If $\bS_0=(1,3,3)$, then  $\bS_0=\set{b_2,\,c_j\succ c_{j+1}\succ c_{j+2}}\cup\bA'$, where
 $\bA'\subset\bA\cup\Pi(b)$. Note that $\bA'$ contains at least two elements from
 $\Pi(b)\cup\set a$. Then $\fD(b_2)=\fD'(b_2)\ge2,\ \fD(a)\ge2$. Hence, $\fD\ge\fC_5$ with
 $\supp\fC_5=\set{a\succ a_2\pr b\succ b_2,\,c\succ c_j\succ c_{j+1}\succ c_{j+2}}$, which is
 impossible. Analogously, $\bS_0=(1,2,5)$ is impossible too,
 which accomplishes the proof of the lemma. 
 \end{proof} 

 Now the implication \emph{1(c)}$\Arr$\emph{1(a)} of Theorem~\ref{main} is easy. Namely,
 let a dimension $\fD$ satisfy \emph{1(c)}. Without loss of geherality, we may suppose $\fD$
 quite sincere. Then either $\nW(\bS)\le3$ or $\fD(0)=1$. In the latter case $\fD$ is obviously
 of finite type. In the former case choose a maximal element $a\in\bS$ as
 stated in Lemma~\ref{32}. Every representation $V$ of dimension $\fD$ without direct summands
 from $\rO(a)$ is isomorphic to $\int_aW$ for a representation $W$ of $\bS^a$. The dimension $\fD'$
 of $W$ is subordinate to $\fD$. Especially it satisfies \emph{1(c)} too; moreover,
 $\fD'(0)=\dim V(a)<\fD(0)$. Thus, using induction by $\fD(0)$, we get that there are finitely
 many nonisomorphic representations of dimension $\fD'$. Since there are finitely many
 subordinate dimension, we obtain the same for the dimension $\fD$.
 \qed

\section{Indecomposable representations}
\label{s4}

 Now we shall prove paragraph \emph{2} of the Main Theorem~\ref{main}. To do it,
 we combine derivations with analogues of some results of \cite{d1,d2} about posets
 of finite type. Namely, we use induction by $|\fD|=\sum_{a\in\hat\bS}\fD(a)$.
 The case $\fD(0)=1$ is obvious. Thus, from now on, we suppose that
 $\fD:\hat\bS\to\mN$ is a dimension of finite type, Theorem~\ref{main} holds
 for every dimension of finite type $\fD'$ of representations of any poset $\bS'$
 such that $|\fD'|<|\fD|$,  and, moreover, $\fD$ is sincere.
 Let $u\in\El_\fD(u)$ be an indecomposable element, $V$ be the corresponding
 representation of $\bS$ and $M$ be the block matrix of the form \eqref{e1}
 describing $u$. Fix an element $a\in\bS$ and denote by $M_a$
 the part of $M$ consisting of the blocks $M(b)$ with $b\pe a$.

  \begin{lemma}\label{41} 
  The columns of the matrix $M_a$ are linear independent.
 \end{lemma} 
 \begin{proof} 
 Obviously, we may suppose that the element $a$ is maximal.
  Consider the part $\ol M$  of $M$ consisting of all blocks $M(b)$ with $b\ne a$. It 
 also describes an object $\ol u\in\El(\sU)$ (certainly, non-sincere) and $|\Dim\ol u|<|\fD|$.
 Hence, Theorem~\ref{main} holds for every indecomposable direct summand $v$ of
 $\ol u$. Especially, the orbit of $v$ is open dense in the space of all objects of the
 same dimension. Let $N$ be the block matrix describing $v$, $\ol\fD=\Dim v$,
 $m=\ol\fD(0)$ and $n=\sum_{b\pr a}\ol\fD(b)$. For every object $v'\in\El_{\ol\fD}(\sU)$,
 denote by $N'$ the corresponding block matrix and by $N_a'$ its part consisting of
 the blocks $N'(b)$ with $b\pr a$.  If $m\le n$, the objects $v'$ such that the rows of
 $N'_a$ are linear independent form an open subset in $\El_{\ol\fD}(\sU)$. Hence,
 $v$ belongs to this subset, i.e. $\mathrm{rk}\,N_a=\ol\fD(0)$.  Then, using automorphisms
 of $u$, one can make zero the part of the matrix $M(a)$ consisting of the rows that
 occur in $v$. Therefore, $v$ is a direct summand of $u$, which is impossible.
 Thus $m>n$. Then the same argument shows that the columns of $v$ are linear
 independent. Since it is so for every direct summand of $\ol u$, it holds for $\ol u$
 too. If, nevertheless, the columns of $M_a$ are linear dependent, then, using
 an automorphism of $u$, one can make a zero column in $M(a)$, which is also
 impossible.
 \end{proof} 

  \begin{corollary}\label{42} 
 For any $a\in\bS$, $\Hom(u,T_a)=0$. Especially,
 neither nonzero endomorphism of $u$ factors through a direct sum of trivial 
 elements; thus $\End u=\End V$, where $V=\rho(u)$ is the corresponding
 representation of $\bS$.
 \end{corollary}

  \begin{lemma}\label{43} 
  Let $\bS$ contain a subset
 $$ 
  
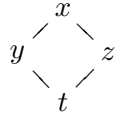
\xymatrix@=1ex{ & x \ar@{-}[dr] \ar@{-}[dl] \\ 
		y \ar@{-}[dr] && z \ar@{-}[dl] \\
		& t	}
 $$ 
 and $V$ be an indecomposable representation of $\bS$ such that its dimension $\fD$
 is of finite type. Then either $\fD(x)=0$ or $\fD(t)=0$.
 \end{lemma} 
 \begin{proof} 
 We may suppose that $\nW(\bS)=3$.
 Again we use the induction by $|\fD|$; for $|\fD|=1$ the claim is trivial. If $a\in\bS$ is
 maximal and $\fD'=\Dim D_aV$, then $\fD'$ is also of finite type and $|\fD'|<|\fD|$.
 If the element $x$ is not maximal, choose $a$ such that $a\succ x$; then $\fD(x)=\fD'(x)$
 and $\fD(t)=\fD'(t)$, so one of them is $0$. Suppose that $x$ is maximal and there is
 another maximal element $a$ such that $t\pr b$. If $p=\set{x,b}\in\Pi(a)$, then $y\pr p$
 and $z\pr p$, hence either $\fD'(t)=\fD(t)=0$ or $\fD'(a)=0$ and $\fD'(p)=0$ for 
 each pair $p=\set{x,b}\in\Pi(a)$, wherefrom $\fD(x)=0$. At last, suppose that $t\nprec a$ for
 any maximal $a\ne x$. If such an element $a$ exists, then $\nW(\Th(a))=2$, hence $\nW(\Th(x))=1$
 and $\bS$ is semidecomposable as $\set x\cup\De'(x)\cup\Th(x)$. Therefore either $\fD(x)=0$
 or $\fD|_{\De'(x)}=0$, especially $\fD(t)=0$.
 \end{proof} 

 The following result is crucial for the proof.

 \begin{lemma}\label{44} 
 A maximal element $a\in\bS$ can be so chosen that  $\End v\iso \End u$
 whenever $u=\int_av$ for an object $v\in\El(\sU^a)$. 
 \end{lemma} 
 \begin{proof} 
  By the matrix description of $\int_av$, the matrix $M$, up to a permutation
 of columns, has the form
 $$ 
  M= \begin{array}{|c|c|c|c|c|l}
  \cline{1-5} &&&&& \\[-2ex] 
	X_0 & X & 0 & Y & Z & \\ 
 \cline{1-5} &&&&& \\[-2ex] 
	0 & 0 & I & I & 0 & , \\
 \cline{1-5}
 \end{array}  
 $$ 
 where $I$ denotes an identity matrix,
 the part $X_0$ is in the matrix $M(a)$ and the remaining part of 
 the first row describes the element $v$. Namely, the part $X$ is in the matrices $M(b),\ b\pr a$;
 the part $Y$ corresponds to the part of $M(b),\ b\in\Th(a)$, arising from the elements $p\in\Pi(a)$
 of $S^a$, such that $b=p''$, while the zero matrix in the first row arises from those $p$ with $b=p'$.
 At last, the part $Z$ arises from the elements $b\in\Th(a)$ considered as the elements of $\bS^a$.
 An endomorphism of $u$ is given by a pair of matrices
 $$ 
  \Phi_0= \begin{pmatrix}
    A & B \\ C & D
 \end{pmatrix}  ,\qquad
  \Phi_1=  \begin{pmatrix}
  S_0 & 0 & 0 & 0 & 0 \\	S_1 & S_2 & S_3 & S_4 & S_5 \\
  0 & 0 & T_{11} & T_{12} & T_{13} \\ 0 & 0 & T_{21} & T_{22} & T_{23} \\
  0 & 0 & T_{31} & T_{32} & T_{33}
 \end{pmatrix}
 $$ 
 such that $\Phi_0 M=M\Phi_1$. Here the diagonal blocks are square, the division of $\Phi_0$
 reflects the horizontal division of $M$, while the division of $\Phi_1$ reflects the vertical division
 of $M$. Note that, by construction, the rows of the matrix $(X_0\ X)$ are linear independent,
 and, by Lemma~\ref{41}, its columns are linear independent too; thus this matrix is invertible.
 It immediately gives that $C=0$. The other equalities for the elements of $\Phi_0M$
 and $M\Phi_1$ are 
 \begin{align*} 
 	AX_0  &= X_0S_0+XS_1 ,\\
	AX	&= XS_2, \\
	B	&= XS_3+YT_{21}+ZT_{31},\\
	B+AY	&= XS_4+YT_{22}+ZT_{32},\\
	AZ	&= XS_5+YT_{23}+ZT_{33}, \\
	D	&= T_{11}+T_{21}=T_{12}+T_{22},\\
	0	&= T_{13}+T_{23}.
 \end{align*} 
 Equivalently, the matrices 
 $$ 
  \Psi_0=A,\quad \Psi_1= \begin{pmatrix}
  	S_2 & S_3-S_4 & S_5 \\  0 & T_{22}-T_{21} & T_{23} \\
	0 & T_{32}-T_{31} & T_{33}
 \end{pmatrix}  
 $$ 
 define an endomorphism of the representation $v$ (then the matrices $S_0,S_1$ can be uniquely
 calculated from the first equality). We have to show that if $\Psi_0=\Psi_1=0$, also $\Phi_0=\Phi_1=0$.
 From the proof of Lemma~\ref{32} we have got to know
 that the element $a$ can be so chosen that every pair from $\Pi(a)$ is of the sort $p=\set{b,c}$,
 where $b$ is a maximal element. If $\set{b,c'}$ is another pair, then neither $b\pe c'$ nor $c\pe b$.
 If an element at some position in the matrix $T_{21}$ is nonzero, it corresponds either to the relation
 $c\pe c'$ or to $b\pr b$; thus the element in the same position of the matrix $T_{22}$ is zero.
 Consequently, the equality $T_{22}-T_{21}=0$ implies that $T_{22}=T_{21}=0$. 
  Lemma~\ref{43} implies that if $t\pr p$, i.e. $t\pr b$ and $t\pr c$, then either $t$-part or
 $p$-part in the element $v$ is empty. It implies, just as above, that if
 $S_3-S_4=0$, then $S_3=S_4=0$, and if $T_{32}-T_{31}=0$,
 then $T_{31}=T_{32}=0$. Thus we get the necessary assertion.
 \end{proof} 

 Now the induction is obvious (just as in \cite{d1,d2}). Namely, if $u$ is an indecomposable
 element such that $\fD$ is of finite type and $u\notin \rO(a)$, then $u=\int D_au$,
 $v=D_au$ is indecomposable, its dimension $\fD'=\Dim v$ is of finite type and
 $|\fD'|<|\fD|$. Therefore, $\End u=\End v=\Mk$, so the stabilizer of the element
 $u$ in the group $\bG_\fD$ is $1$-dimensional. Recall that there always is an open
 subset $U\sbe\El_{\fD}(\sU)$ such that the stabilizers of the elements of $U$ are of
 minimal dimension (see e.g. \cite{sh}). Thus $u\in U$. Since any orbit is open in its
 closure, we get that the orbit of $u$ is open (hence dense), so its dimension, which is
 $\dim\bG_\fD-1$, equals $\dim\El_\fD(\sU)$. Therefore
 $\rQ_\bS(\fD)=\dim\bG_\fD-\dim\El_{\fD}(\sU)=1$.
 Moreover, if $u'$ is another indecomposable element of the same dimension $\fD$,
 its orbit is also open, hence coincide with that of $u$, i.e. $u$ is the unique indecomposable
 element of this dimension. 

 On the other hand, if a dimension is of finite type and $\rQ_\bS(\fD)=1$, the number of orbits
 of the group $\bG_\fD$ in $\El_\fD(\sU)$ is finite. Therefore, there is an open orbit and its
 dimension equals $\dim\El_\fD(\sU)$. If $u$ is an element of this orbit,
 $\dim\End u=\dim\bG_\fD-\dim_\fD(\sU)=\rQ_\bS(\fD)=1$, so $\End u=\Mk$ and $u$ is
 indecomposable. It accomplishes the proof of the Main Theorem.

\end{document}